\def\Bbb#1{{\fam\msbfam\relax#1}}
\font\fivemsb=msbm5
\font\sevenmsb=msbm7
\font\tenmsb=msbm10
\newfam\msbfam
\textfont\msbfam=\tenmsb
\scriptfont\msbfam\sevenmsb
\scriptscriptfont\msbfam\fivemsb

\def\spc{{\Bbb C}}
\def\spr{{\Bbb R}}

\def\wt{\widetilde}
\def\wh{\widehat}
\def\a{\alpha}
\def\b{\beta}

\def\g{\gamma}
\def\l{\lambda}

\def\r{\rho}
\def\s{\sigma}
\def\p{\pi}

\def\t{\theta}
\def\x{\xi}
\def\y{\eta}

\def\ce{{\cal E}}
\def\cg{{\cal G}}
\def\ck{{\cal K}}
\def\cl{{\cal L}}

\def\la{\langle}
\def\ra{\rangle}

\documentstyle[verbatim]{amsart}

\theoremstyle{plain}
\newtheorem{theorem}{Theorem}[section]

\newtheorem{lemma}[theorem]{Lemma}
\newtheorem{proposition}[theorem]{Proposition}

\theoremstyle{definition}
\newtheorem{definition}[theorem]{Definition}
\newtheorem{example}[theorem]{Example}

\numberwithin{equation}{section}

\renewcommand{\appendixname}{\sc Appendix}

\begin{document}
\title[Log-concavity, log-convexity, growth order]%
   {Roles of Log-concavity, log-convexity, \\and growth order
    in white noise analysis}

\maketitle
\begin{center}
{\sc Nobuhiro Asai}\\
{\it 
International Institute for Advanced Studies\\
Kyoto, 619-0225, JAPAN.}
\end{center}
\smallskip
\begin{center}
{\sc Izumi Kubo}\\
{\it Department of Mathematics 
\\Graduate School of Science\\
Hiroshima University \\
Higashi-Hiroshima, 739-8526, JAPAN}
\end{center}
\begin{center}
and
\end{center}
\begin{center}
{\sc Hui-Hsiung Kuo}\\
{\it Department of Mathematics\\
Louisiana State University\\
Baton Rouge, LA 70803, USA}
\end{center}

\bigskip
\medskip

\begin{abstract}
In this paper we will develop a systematic method 
to answer the questions 
$(Q1)(Q2)(Q3)(Q4)$ (stated in Section 1)
with complete generality.  As a result,
we can solve the difficulties $(D1)(D2)$ 
(discussed in Section 1)
without uncertainty.  For these purposes 
we will introduce certain classes of growth
functions $u$ and apply the Legendre transform 
to obtain a sequence which leads to
the weight sequence $\{\a(n)\}$ first studied by 
Cochran et al. \cite{cks}.  The notion of 
(nearly) equivalent functions, (nearly) equivalent 
sequences and dual Legendre functions will be defined  
in a very natural way.  An application to the growth order of 
holomorphic functions on $\ce_c$ will also be discussed.
\end{abstract}

\section{Introduction} \label{sec:1}

Let $\ce$ be a real nuclear space with topology given by a
sequence of inner product norms $\{|\cdot|_{p}
\}_{p=0}^{\infty}$. Let $\ce_{p}$ be the completion of $\ce$
with respect to the norms $|\cdot|_{p}$. 
We will assume the following conditions:
\begin{itemize}
\item[(a)] There exists a constant $0<\r<1$ such that
$|\cdot|_0\leq \r|\cdot|_{1}\leq \cdots\leq \r^{p}|\cdot|_{p}
\leq \cdots$.
\item[(b)] For any $p\geq 0$, there exists some $q\geq p$ such
that the inclusion mapping $i_{q, p}: \ce_{q} \to \ce_{p}$ is
a Hilbert-Schmidt operator.
\end{itemize}
We denote the complexification of $\ce$ by $\ce_c$. 

The generalized and test functions in 
$[\ce]_{\a}^{*}$ and $[\ce]_{\a}$ (CKS-space for short) 
are characterized in
terms of their $S$-transforms $F$
under a 
very general setting by \cite{cks} and \cite{akk2}, respectively. 
There are two conditions on
$F$. The first one is the analyticity.
The second one is the growth condition. 
The exponential generating functions
\begin{equation}
	G_\a(r)=\sum_{n=0}^{\infty}{\a(n)\over n!}r^n, \qquad
	G_{1/\a}(r)=\sum_{n=0}^{\infty}{1\over n!\a(n)}r^n
\end{equation}
are adopted as the growth functions 
for generalized functions \cite{cks}
and test functions \cite{akk2}, respectively. 
Recently, Asai et al. have shown in \cite{akk4} that  
the conditions $(A1)(A2)(B2)
(\widetilde{B}2)(C2)(\widetilde{C}2)$ (See Appendix)
are minimal assumptions on $\a(n)$ not only for the construction of 
CKS-space, but also for white noise operator theory on it 
(See also \cite{akk3}\cite{kks}\cite{ob99}).

On the other hand, the recent paper \cite{ghor} by Gannoun et al.
is closely connected  with the paper \cite{cks} by Cochran et al. 
and the series of papers 
\cite{akk1}\cite{akk2}\cite{akk3}\cite{akk4}\cite{akk5}.
They study the spaces of holomorphic functions 
$\cal G_{\t^{*}}$ on $\ce_c$ and $\cal F_{\t}$
on the dual space $\ce'_{c}$ 
considering functions 
$\exp (\t(r)^{*})$ and $\exp(\t(r))$
as growth functions where $\t^{*}$ is the dual function
of Young function $\t$.  They remark briefly the relationship 
between $\cal G_{\t^{*}}$ and $[\ce]_{\a}^{*}$ by 
taking $\t(r)^{*}=\log G_{\a}(r^2)$ and 
$\a(n):=n!\t_n=\inf_{r>0}{\exp(\t(r)^{*})\over r^n}$.

However, the papers \cite{cks} and \cite{ghor} 
do not provide  a general method to solve
the following delicate problems:
\begin{itemize}
\item[(D1)]
Technical difficulties mentioned 
after $(D2)$ cannot be canceled out by systematic way.
Moreover, 
Young function $\t(r)$ cannot be obtained explicitly from
$\t(r)^{*}=\log G_{\a}(r^2)$ even for the important cases,  
$\a(n)=(n!)^{\b}$ and Bell's numbers $b_k(n)$ given in \cite{cks}. 
\item[(D2)]
It is not checked whether $n!\t_n$ satisfies
$(A2)(B2)(\widetilde{B}2)(C2)(\widetilde{C}2)$.  
Since precise estimates are required in general,
this is not obvious problem at all.  
Consult papers \cite{akk1}\cite{kks} for the case of 
Bell's numbers.
\end{itemize}

Now we shall justify our claim $(D1)$ as follows.
In \cite{cks} the following growth condition is used
for generalized functions:
\begin{itemize}
\item[$\bullet$] There exist constants $K, a, p\geq 0$ such
that
\begin{equation} \label{eq:1-10}
	|F(\x)| \leq K G_{\a}\big(a|\x|_{p}^{2}\big)^{1/2},
        \qquad \x\in\ce_{c}.
\end{equation}

\end{itemize}
On the other hand, in \cite{akk2} the
following growth condition is used for test functions:
\begin{itemize}
\item[$\bullet$] For any constants $a, p\geq 0$, there exists
a constant $K\geq 0$ such that
\begin{equation} \label{eq:1-11}
|F(\x)| \leq K G_{1/\a}\big(a|\x|_{-p}^{2}\big)^{1/2},
\qquad \x\in\ce_{c}.
\end{equation}
\end{itemize}

In the case of Kondratiev-Streit space, $G_{\a}(r)$ and
$G_{1/\a}(r)$ are given by
\begin{equation} \label{eq:e}
G_{\a}(r) = \sum_{n=0}^{\infty} {1 \over (n!)^{1-\b}}\,
  r^{n}, \quad
G_{1/\a}(r) = \sum_{n=0}^{\infty} {1 \over (n!)^{1+\b}}\,
  r^{n}.
\end{equation}
These series cannot be summed up in closed forms unless
$\b=0$ (the case of Hida-Kubo-Takenaka space). 
Fortunately, we have the estimates:
\begin{align}
\exp\Big[(1-\b)\,r^{{1\over 1-\b}}\Big] & \leq G_{\a}(r)
  	\leq 2^{\b}\exp\Big[(1-\b)\, 2^{{\b\over 1-\b}}\,
  	r^{{1\over 1-\b}}\Big]. \label{eq:c} \\
	2^{-\b}\exp\Big[(1+\b)\, 2^{-{\b\over 1+\b}}\,
  	r^{{1\over 1+\b}}\Big] & \leq G_{1/\a}(r)  \leq
\exp\Big[(1+\b)\,r^{{1\over 1+\b}}\Big]. \label{eq:d}
\end{align}
Therefore we can substitute the growth functions $G_{\a}$ and
$G_{1/\a}$ in Equations (\ref{eq:1-10}) and (\ref{eq:1-11})
by the following functions $\wt G_{\a}$ and $\wt G_{1/\a}$,
respectively,
\begin{equation} \label{eq:1-13}
\wt G_{\a}(r) = \exp\Big[(1-\b) r^{1\over 1-\b}\Big], \qquad
\wt G_{1/\a}(r) = \exp\Big[(1+\b) r^{1\over 1+\b}\Big].
\end{equation}
These are the growth functions used in \cite{ks92}
\cite{ks93} (see \cite{ps}
and \cite{kps} for $\b=0$, respectively.)

Now, for the case of CKS-space,
the growth functions
$G_{\a}$ and $G_{1/\a}$ in Equations (\ref{eq:1-10}) and
(\ref{eq:1-11}) are not practical to use since in general
we cannot find closed forms for the sums of the infinite series.
For example, for
the case of $\a(n)=b_{k}(n)$,
even though $G_{\a}(r)=\exp_{k}(r)$, the $k$-th
iterated exponential function, we simply do not have a
closed form for the corresponding $G_{1/\a}$. 

The main purpose of this paper is to
consider the following questions:
\begin{itemize}
\item[(Q1)] 
In Equations (\ref{eq:1-10}) and
(\ref{eq:1-11}) can we replace the growth functions $G_{\a}$
and $G_{1/\a}$ by elementary functions $U$ and $u$,
respectively?
\item[(Q2)]
How to find $U$ and $u$ from $G_{\a}$ and
$G_{1/\a}$? In particular, how to find $\wt G_{\a}(r)$ and
$\wt G_{1/\a}(r)$ in Equation (\ref{eq:1-13}) from $G_{\a}$
and $G_{1/\a}$ in Equation (\ref{eq:e}) and without
appealing to Equations (\ref{eq:c}) and (\ref{eq:d})?
This is related to $(D1)$.
\item[(Q3)]
Are there any other general criteria to check 
$(A2)(B2)(\widetilde{B}2)(C2)(\widetilde{C}2)$
for a given growth function $U$ or $u$ 
without technical estimates?
This is connected with $(D2)$.
\item[(Q4)] 
For $(Q1)\sim (Q3)$ what kinds of conditions do we have to 
impose on functions $U$ and $u$?
\end{itemize}

The Legendre transform is the key tool to solve 
the above four questions.
The answers to $(Q3)$ can be found in 
Theorems \ref{lem:3-2}, \ref{lem:3-4} 
and \ref{lem:a-1}.  
For $(Q1)(Q2)$, Theorems \ref{thm:3-1}, 
\ref{thm:4-2} and \ref{thm:4-3} will play 
fundamental roles.  
About $(Q4)$, please refer to Definition \ref{def:2-1},
Equations 
\eqref{eq:3-3} and \eqref{eq:4-5} for quick reference.   

The present paper is organized as follows. 
In Section 2, several kinds of log-convex functions
will be prepared for the Legendre transform.  
In Section 3, we examine 
properties of the Legendre transform for various log-convex
functions.   In addition, we will introduce the notion of 
equivalent functions and sequences motivated from 
Theorem \ref{thm:3-1}.
In Section 4, we discuss the dual Legendre 
functions, nearly equivalent functions and sequences.
In Section 5, an application to the growth order of 
holomorphic functions on $\ce_c$ will be discussed 
under quite general assumptions.  In paticular, 
Theorem \ref{thm:3-1} 
will be useful to prove the topological 
isomorphism between spaces ${\cal G}_u$ and ${\cal K}_u$. 
Further related works with the present paper 
can be found in \cite{akk3} and \cite{akk5}.

\section{Log-, (log, exp)-, and (log, $x^{k}$)-convex
functions} \label{sec:2}

In this section, we shall consider three kinds of convexity 
for later use. Before giving their definitions,
let us start with the following, which stems
from the proof of Theorem 4.3 in \cite{cks}
and is connected with $(Q1)(Q2)$.   
It explains our viewpoint.

Let $u(x) =
\sum_{n=0}^{\infty} u_{n} x^{n}$ be an entire function with
$u_{n}>0$ and the sequence $\{u_{n}\}$ being log-concave
(See $(C2)(C3)$ in Appendix.),
i.e.,
\begin{equation} \label{eq:h}
u_{n} u_{n+2} \leq u_{n+1}^{2},\qquad\forall n\geq 0.
\end{equation}
It is shown in the proof of Theorem 4.3 in \cite{cks} that
\begin{equation}
1 \leq {1\over u_{n}} \inf_{r>0} {u(r) \over r^{n}}
\leq (n+1)\,e.  \notag
\end{equation}
But $n+1\leq 2^{n}$. Hence
\begin{equation} \label{eq:cks}
u_{n} \leq \inf_{r>0} {u(r) \over r^{n}}
  \leq e\, 2^{n}u_{n}.
\end{equation}
These inequalities imply that
\begin{equation} \label{eq:k}
u(x) \leq \sum_{n=0}^{\infty} \left(\inf_{r>0} {u(r) \over
r^{n}}\right) x^{n} \leq e\, u(2x).
\end{equation}
Thus if we can regard the last series as the function
$G_{1/\a}$, then the function $u$ can be used as the growth
function for test functions. But to get a satisfactory
answer, we need to describe $u$ by function properties
instead of expressing it as an infinite series. This leads to
a general question:
\begin{itemize}
\item[(Q)] What functions $U$ and $u$ can serve as growth
functions in Equations (\ref{eq:1-10}) and (\ref{eq:1-11})
for generalized and test functions, respectively?
\end{itemize}

Now we are in a position to define the three kinds of 
log-convexity of functions as follow.
\begin{definition}\label{def:2-1}
Let $u$ be a positive continuous function on $[0, \infty)$.
\begin{itemize}
\item[(a)] The function $u$ is called {\em log-convex} if
$\log u$ is convex on $[0, \infty)$.
\item[(b)] The function $u$ is called {\em (log, exp)-convex}
if $\log u(e^{x})$ is convex on $\spr$.
\item[(c)] The function $u$ is called
{\em (log, $x^{k}$)-convex} if $\log u(x^{k})$ is convex on
$[0, \infty)$. Here $k$ is a positive real number.
\end{itemize}
\end{definition}

The concept of (log, exp)-convexity will play an important
role in this paper. The (log, exp)-convex functions are
the appropriate functions to replace those given by infinite
series $u(r)=\sum_{n=0}^{\infty} u_{n}r^{n}$ as mentioned 
above. In fact, it has been shown in \cite{kubo}
that if $u(r)=\sum_{n=0}^{\infty} u_{n}r^{n}$ is an entire
function with $u_{n}\geq 0$ and $u(r)>0$ for all $r\geq 0$,
then the function $u$ is a (log, exp)-convex function. This
fact also follows from Equation (\ref{eq:2-1}) below and the
equality
\begin{align}
& u(r)u''(r)-u'(r)^{2} + {1\over r}u(r)u'(r) \notag \\
& = u_{0}u_{1}{1\over r} + \sum_{n=0}^{\infty}\left(
  \sum_{j=0}^{[(n+2)/2]}
   (n+2-2j)^{2} u_{j}u_{n+2-j}\right) r^{n}. \notag
\end{align}
However, a (log, exp)-convex function $u$ may not be given
by an entire function with positive coefficients in the
series expansion. For instance, $u(r)=\exp[r^{2}-r^{3}+
r^{4}]$ is such a function.

\begin{example} \label{exa:2-1}
Suppose $u$ is a positive $C^{2}$-function on $[0, \infty)$.
It is easy to check by direct calculations the following
assertions:
\begin{itemize}
\item[(1)] $u$ is log-convex if and only if
\begin{equation}
u(r)u''(r)-u'(r)^{2} \geq 0, \qquad \forall r > 0. \notag
\end{equation}
\item[(2)] $u$ is (log, exp)-convex if and only if
\begin{equation} \label{eq:2-1}
u(r)u''(r)-u'(r)^{2} + {1\over r}u(r)u'(r) \geq 0,
\qquad \forall r > 0.
\end{equation}
\item[(3)] $u$ is (log, $x^{k}$)-convex if and only if
\begin{equation}
u(r)u''(r) -u'(r)^{2}+{k-1 \over k r} u(r)u'(r)
\geq 0,  \qquad \forall r > 0. \notag
\end{equation}
\end{itemize}
\end{example}

Observe from this example that if $u$ is log-convex then
it is convex. When $u$ is an increasing function,
we have the implications: (i) log-convex $\!\!\implies\!\!$
(log, $x^{k}$)-convex for any $k\geq 1$,
(ii) (log, $x^{k}$)-convex for some $k>0$ $\!\!\implies\!\!$
(log, exp)-convex. In fact, these implications are true
in general.

\begin{proposition} \label{pro:2-1}
Let $u$ be a positive continuous function on $[0, \infty)$.
\begin{itemize}
\item[(1)] If $u$ is log-convex, then it is convex.
\end{itemize}
When $u$ is also increasing, we have the assertions:
\begin{itemize}
\item[(2)] If $u$ is log-convex, then it is
(log, $x^{k}$)-convex for any $k\geq 1$.
\item[(3)] If $u$ is (log, $x^{k}$)-convex for some $k>0$,
then it is (log, exp)-convex.
\end{itemize}
\end{proposition}

\begin{pf}
To prove assertion (1), let $r, s\geq 0$ and $\l\in [0, 1]$.
Since by assumption $u$ is log-convex, we have
\begin{equation} \label{eq:f}
\log u(\l r + (1-\l)s) \leq \l \log u(r) + (1-\l)\log u(s).
\end{equation}
But $e^{x}$ is a convex function. Hence
\begin{align}
e^{\l \log u(r)+(1-\l)\log u(s)}
  & \leq \l e^{\log u(r)} + (1-\l) e^{\log u(s)} \notag \\
  & = \l u(r) + (1-\l) u(s). \label{eq:g}
\end{align}
Take the exponential in Equation (\ref{eq:f}) and then use
Equation (\ref{eq:g}) to get
\begin{equation}
u(\l r + (1-\l)s) \leq \l u(r) + (1-\l)u(s).  \notag
\end{equation}
Hence $u$ is convex and we have proved assertion (1). For
assertion (2) we use the fact that $x^{k}$ is convex for any
$k\geq 1$ and the assumption that $u$ is increasing to get
\begin{equation} \label{eq:2-2}
u\big((\l r+(1-\l)s)^{k}\big) \leq u\big(\l r^{k}+
(1-\l)s^{k}\big).
\end{equation}
Suppose $u$ is log-convex. Then
\begin{equation} \label{eq:2-3}
\log u\big(\l r^{k}+(1-\l)s^{k}\big) \leq \l \log u(r^{k})
  + (1-\l) \log u(s^{k}).
\end{equation}
Upon taking logarithm in Equation (\ref{eq:2-2}) and then use
Equation (\ref{eq:2-3}) we obtain
\begin{equation}
\log u\big((\l r+(1-\l)s)^{k}\big) \leq \l \log u(r^{k})
  + (1-\l) \log u(s^{k}). \notag
\end{equation}
This shows that $\log u(x^{k})$ is convex, i.e., $u$ is
(log, $x^{k}$)-convex and so assertion (2) is proved.
For the third assertion, note that $e^{x}$ is convex and
the function $u(x^{k})$ is increasing. Hence we have
\begin{equation} \label{eq:2-4}
u\Big(\big(e^{\l r/k +(1-\l) s/k}\big)^{k}\Big)
\leq u\Big(\big(\l e^{r/k} + (1-\l) e^{s/k}\big)^{k}\Big).
\end{equation}
But by assumption $u$ is (log, $x^{k}$)-convex. Hence
\begin{equation} \label{eq:2-5}
\log u\Big(\big(\l e^{r/k} + (1-\l) e^{s/k}\big)^{k}\Big)
\leq \l \log u\big((e^{r/k})^{k}\big) +
  (1-\l) \log u\big((e^{s/k})^{k}\big).
\end{equation}
Upon taking logarithm in Equation (\ref{eq:2-4}) and then use
Equation (\ref{eq:2-5}) we obtain
\begin{equation}
\log u\big(e^{\l r + (1-\l)s}\big) \leq \l \log u(e^{r})
+ (1-\l)\log u(e^{s}). \notag
\end{equation}
Hence the function $\log u(e^{x})$ is convex, i.e., $u$ is
(log, exp)-convex and so assertion (3) is proved.
\end{pf}

\begin{lemma} \label{lem:b-1}
Let $f$ be a convex function on $\spr$ such that
$\lim_{x\to -\infty} f(x)$ exists. Then the function $f$
is an increasing function.
\end{lemma}

\begin{pf}
Let $x<x_{1}<x_{2}$ with $x_{1}$ and $x_{2}$ being fixed.
Then
\begin{equation}
x_{1} = {x_{2}-x_{1} \over x_{2}-x} x +
  {x_{1}-x \over x_{2}-x} x_{2}.  \notag
\end{equation}
Since the function $f$ is convex,
\begin{equation}
f(x_{1}) \leq {x_{2}-x_{1} \over x_{2}-x} f(x) +
  {x_{1}-x \over x_{2}-x} f(x_{2}). \notag
\end{equation}
Letting $x\to -\infty$, we get $f(x_{1}) \leq f(x_{2})$
and so the function $f$ is increasing.
\end{pf}

\begin{lemma} \label{lem:b-2}
If a positive continuous function $u$ on $[0, \infty)$ is
(log, exp)-convex, then it is an increasing function.
\end{lemma}

\begin{pf}
Let $f(x)=\log u(e^{x})$. By Lemma \ref{lem:b-1}, $f$ is
increasing. It follows that the function $u$ is also
increasing.
\end{pf}

We want to point out that $u(r)$ being defined at $r=0$
is crucial for Lemma \ref{lem:b-2}. For example, let
$u(r)=\exp\big[(\log r)^{2}-2\log r\big]$. Obviously, the
function $\log u(e^{x})$ is convex on $\spr$. But $u$ is
not an increasing function on $(0, \infty)$.

\section{Legendre transform} \label{sec:3}

For the characterization theorems in the paper by Cochran
et a.~\cite{cks} the following condition is imposed
\begin{equation}
\limsup_{n\to\infty} \left({n! \over \a(n)}\,\inf_{r>0}
{u(r) \over r^{n}}\right)^{1/n} < \infty. \notag
\end{equation}
This condition leads to the consideration of the sequence
\begin{equation} \label{eq:3-1}
\inf_{r>0} {u(r) \over r^{n}}, \qquad n=0, 1, 2, \ldots.
\end{equation}
Moreover, as we pointed out in Section \ref{sec:2} that if
$u(x)=\sum_{n=0}^{\infty} u_{n} x^{n}$ is an entire function
with $u_{n}>0$ and $\{u_{n}\}$ being log-concave, then
Equation (\ref{eq:k}) holds, i.e., the functions $u$ is
``equivalent'' to the following function
\begin{equation} \label{eq:3-2}
\sum_{n=0}^{\infty} \left(\inf_{r>0} {u(r) \over
r^{n}}\right) x^{n}.
\end{equation}

The above discussion raises three questions: (i) What
function $u$ can we define the sequence in Equation
(\ref{eq:3-1})? (ii) What is the new function in Equation
(\ref{eq:3-2})? (iii) Is this new function ``equivalent'' to
the function $u$? In this section we will give answers to
these questions.

\smallskip
\noindent
{\bf Notation.} Let $C_{+, \log}$ denote the set of all
positive continuous functions $u$ on $[0, \infty)$ satisfying
the condition
\begin{equation} \label{eq:3-3}
\lim_{r\to \infty} {\log u(r) \over \log r}=\infty.
\end{equation}

\smallskip
Observe that the condition in Equation (\ref{eq:3-3}) means
that $u$ grows faster than all polynomials. The set
$C_{+, \log}$ includes all entire functions $u(r)=
\sum_{n=0}^{\infty} u_{n} r^{n}$ with $u_{n}\geq 0$ for all
$n$ and $u_{n}>0$ for infinitely many $n$'s. If $u$ is
a function in the set $C_{+, \log}$, then we can define the
sequence in Equation (\ref{eq:3-1}). In fact, we
will define Equation (\ref{eq:3-1}) as a function on
$[0, \infty)$.

\begin{definition}
The {\em Legendre transform} $\,\ell_{u}$ of a function
$u\in C_{+, \log}$ is defined to be the function
\begin{equation} \label{eq:3-4}
\ell_{u}(t) = \inf_{r>0} {u(r) \over r^{t}},
  \qquad t\in [0, \infty).
\end{equation}
\end{definition}

The next lemma is immediate from the definition of Legendre
transform.

\begin{lemma} \label{lem:3-1}
{\em (1)} For $u\in C_{+, \log}$ and $a>0$, let $\t_{a}u$
be the function $\t_{a}u(r)=u(ar)$. Then $\t_{a}u\in
C_{+, \log}$ and $\ell_{\t_{a}u}(t) = a^{t} \ell_{u}(t)$
for all $t\geq 0$.

\noindent
{\em (2)} Suppose $u, v\in C_{+, \log}$ and $u(r)\leq v(r)$
for all $r\geq 0$. Then $\ell_{u}(t)\leq \ell_{v}(t)$ for
all $t\geq 0$.
\end{lemma}

\begin{definition}
A positive function $f$ on $[0, \infty)$ is called
{\em log-concave} if $\log f$ is a concave function,
or equivalently, for any $t_{1}, t_{2}\geq 0$ and $0\leq
\l \leq 1$,
we have
\begin{equation} \label{eq:3-5}
f(\l t_{1} + (1-\l)t_{2}) \geq f(t_{1})^{\l}\,f(t_{2})^{1-\l}.
\end{equation}
\end{definition}

Put $t_{1}=n, t_{2}=n+2$, and $\l=1/2$ in Equation
(\ref{eq:3-5}) to get
\begin{equation}
f(n) f(n+2) \leq f(n+1)^{2}, \qquad\forall n\geq 0.  \notag
\end{equation}
This shows that if a positive function $f$ on $[0, \infty)$
is log-concave, then the sequence $\{f(n)\}$ is log-concave
(see Equation (\ref{eq:h}).)

\begin{theorem} \label{lem:3-2}
The Legendre transform $\ell_{u}$ of a function $u\in
C_{+, \log}$ is log-concave. {\em (}Hence the function
$\ell_{u}(t)$ is continuous and the sequence
$\{\ell_{u}(n)\}_{n=0}^{\infty}$ is log-concave.{\em )}
\end{theorem}

\begin{pf}
For any $t_{1}, t_{2}\geq 0$ and $0\leq \l \leq 1$, we have
\begin{align}
\ell_{u}(\l t_{1} + (1-\l)t_{2})
& = \inf_{r>0} {u(r) \over r^{\l t_{1} + (1-\l)t_{2}}}
           \notag \\
& = \inf_{r>0} {u(r)^{\l} \over r^{\l t_{1}}}\,
     {u(r)^{1-\l} \over r^{(1-\l) t_{2}}}  \notag \\
& \geq \left(\inf_{r>0} {u(r) \over r^{t_{1}}}\right)^{\l}
    \left(\inf_{r>0} {u(r) \over r^{t_{2}}}\right)^{1-\l}
    \notag  \\
& = \ell_{u}(t_{1})^{\l}\,\ell_{u}(t_{2})^{1-\l}. \notag
\end{align}
Hence by Equation (\ref{eq:3-5}) the function $\ell_{u}$
is log-concave.
\end{pf}

Now, we consider those functions in $C_{+, \log}$ which are
(log, exp)-convex. Let $u$ be such a function. Then the
left-hand derivative $u_{-}'(r)$ and the right-hand
derivative $u_{+}'(r)$ exist. For convenience, define
\begin{equation}
\tau_{-}(r) = {ru_{-}'(r) \over u(r)}, \qquad
\tau_{+}(r) = {ru_{+}'(r) \over u(r)}.  \notag
\end{equation}
Both $\tau_{-}$ and $\tau_{+}$ are increasing functions.
Since $u$ is increasing by Lemma \ref{lem:b-2}, we have
$0\leq \tau_{-}(r)\leq \tau_{+}(r)$ for all $r\geq 0$.
Moreover, the condition in Equation (\ref{eq:3-3}) implies
that $\tau_{-}(r)\to\infty$ as $r\to\infty$. Note that
for any $r\geq 0$ and $t\in [\tau_{-}(r), \tau_{+}(r)]$,
\begin{equation}
{u(s) \over s^{t}} \geq {u(r) \over r^{t}}, \qquad
\forall s > 0.  \notag
\end{equation}
Hence $\inf_{s>0} u(s)/s^{t}=u(r)/r^{t}$ and so
\begin{equation}
\ell_{u}(t) = {u(r) \over r^{t}}, \qquad
    \text{~for any~} t\in [\tau_{-}(r), \tau_{+}(r)].
  \notag
\end{equation}
In particular, let $\r(t)$ be a solution of the
equation $\tau_{-}(r)=t$, i.e., $\tau_{-}(\r(t))=t$.
Then we have
\begin{equation}
\ell_{u}(t) = {u(\r(t)) \over \r(t)^{t}}. \notag
\end{equation}

We sum up the above discussion in the next lemma.

\begin{lemma} \label{lem:3-3}
Let $u\in C_{+, \log}$ be (log, exp)-convex. Then
\begin{equation} \label{eq:3-6}
\ell_{u}(t) = {u(r) \over r^{t}}, \qquad
   \text{~for any~} t\in [\tau_{-}(r), \tau_{+}(r)],
\end{equation}
where $\tau_{-}(r)=ru_{-}'(r)/u(r)$ and $\tau_{+}(r)=
ru_{+}'(r)/u(r)$. In particular, let $\r(t)$ be a solution
of the equation $\tau_{-}(r)=t$, i.e., $\tau_{-}(\r(t))=t$.
Then
\begin{equation} \label{eq:3-7}
\ell_{u}(t) = {u(\r(t)) \over \r(t)^{t}}.
\end{equation}
\end{lemma}

\begin{theorem} \label{lem:3-4}
Let $u\in C_{+, \log}$ be (log, exp)-convex. Then
its Legendre transform $\ell_{u}(t)$ is decreasing for large
$t$ and $\lim_{t\to\infty} \ell_{u}(t)^{1/t}=0$.
\end{theorem}

\begin{pf}
Let $s\geq t$ be fixed. Use Lemma \ref{lem:3-3} to get
\begin{equation}
\ell_{u}(t)={u(\r(t)) \over \r(t)^{t}}=\r(t)^{s-t}
{u(\r(t)) \over \r(t)^{s}}.  \notag
\end{equation}
Then by the definition of the Legendre transform,
\begin{equation} \label{eq:3-8}
\ell_{u}(t) \geq \r(t)^{s-t} \ell_{u}(s).
\end{equation}
Recall that $\tau_{-}(r)$ increases to $\infty$
monotonically. Hence $\r(t)$ also increases to $\infty$
monotonically. Choose $t_{0}$ such that $\r(t)>1$ for all
$t\geq t_{0}$. Then it follows from Equation (\ref{eq:3-8})
that
\begin{equation}
\ell_{u}(t) \geq \ell_{u}(s), \qquad \forall
  s\geq t\geq t_{0}.  \notag
\end{equation}
Hence $\ell_{u}(t)$ is decreasing for large $t$. Moreover,
from Equation (\ref{eq:3-8}) we have
\begin{equation}
\ell_{u}(s) \leq \r(t)^{t-s} \ell_{u}(t), \qquad \forall
s\geq t\geq t_{0}.  \notag
\end{equation}
Therefore,
\begin{equation}
\ell_{u}(s)^{1/s} \leq \r(t)^{t/s -1}\,\ell_{u}(t)^{1/s},
  \qquad \forall s\geq t\geq t_{0}.  \notag
\end{equation}
Hold $t$ fixed and let $s\to\infty$ to get
\begin{equation}
\limsup_{s\to\infty} \ell_{u}(s)^{1/s} \leq \r(t)^{-1},
\qquad \forall t\geq t_{0}. \notag
\end{equation}
But $\r(t)\to\infty$ as $t\to\infty$. Hence we can conclude
that $\lim_{s\to\infty} \ell_{u}(s)^{1/s}=0$.
\end{pf}

\begin{lemma} \label{lem:3-6}
Let $u\in C_{+, \log}$ be (log, exp)-convex. Then
\begin{equation}
u(r) = \sup_{t\geq 0} \ell_{u}(t) r^{t}, \qquad
  \forall r \geq 0. \notag
\end{equation}
\end{lemma}

\begin{pf}
From the definition of the Legendre transform in Equation
(\ref{eq:3-4}) we have
\begin{equation} \label{eq:3-11}
\ell_{u}(t) r^{t} \leq u(r), \qquad \forall r\geq 0.
\end{equation}
On the other hand, for any fixed $r\geq 0$, we can
choose $t\in [\tau_{-}(r), \tau_{+}(r)]$ in Lemma
\ref{lem:3-3} to get
\begin{equation}  \label{eq:3-12}
\ell_{u}(t) r^{t} = u(r).
\end{equation}
Equations (\ref{eq:3-11}) and (\ref{eq:3-12}) imply that
$u(r)=\sup_{t\geq 0} \ell_{u}(t) r^{t}$.
\end{pf}

The next lemma follows immediately from Lemmas
\ref{lem:3-1} (2) and \ref{lem:3-6}.

\begin{lemma} \label{lem:3-7}
Let $u, v \in C_{+, \log}$ be (log, exp)-convex. Then
\begin{itemize}
\item[(1)] $u=v$ if and only if $\ell_{u} = \ell_{v}$.
\item[(2)] $\ell_{u}(t) \leq \ell_{v}(t)$ for all
$t\geq 0$ if and only if $u(r)\leq v(r)$ for all $r\geq 0$.
\end{itemize}
\end{lemma}

Now, we consider functions $u$ in $C_{+, \log}$ which are
(log, $x^{k}$)-convex.

\begin{lemma} \label{lem:3-5}
Let $u\in C_{+, \log}$ and $k>0$. Then $u$ is
(log, $x^{k}$)-convex if and only if $\ell_{u}(t)t^{kt}$ is
log-convex.
\end{lemma}

\begin{pf}
Let $t_{1}, t_{2}\geq 0$ and $0\leq \l\leq 1$. Then
\begin{equation}
\ell_{u}(t_{1})^{\l} \ell_{u}(t_{2})^{1-\l}
= \inf_{r, s>0} \left({u(r) \over r^{t_{1}}}\right)^{\l}
   \left({u(s) \over s^{t_{2}}}\right)^{1-\l}
= \inf_{x, y>0} {u(x^{k})^{\l}\,u(y^{k})^{1-\l}
     \over x^{\l kt_{1}}\,y^{(1-\l)kt_{2}}}.  \notag
\end{equation}
Suppose $u$ is (log, $x^{k}$)-convex. Then
\begin{equation}
u\big((\l x + (1-\l)y)^{k}\big) \leq u(x^{k})^{\l}\,
   u(y^{k})^{1-\l}.  \notag
\end{equation}
Therefore,
\begin{equation}
\ell_{u}(t_{1})^{\l} \ell_{u}(t_{2})^{1-\l} \geq
\inf_{x, y>0} {u\big((\l x + (1-\l)y)^{k}\big)
     \over x^{\l kt_{1}}\,y^{(1-\l)kt_{2}}}.  \notag
\end{equation}
Make a change of variables $z=\l x + (1-\l)y$ to get
\begin{align}
& \ell_{u}(t_{1})^{\l} \ell_{u}(t_{2})^{1-\l} \notag \\
& \geq (1-\l)^{(1-\l)kt_{2}} \inf_{z>0} u(z^{k})
\inf_{0<\l x <z} \,{1\over x^{\l kt_{1}}
  (z-\l x)^{(1-\l)kt_{2}}}.  \label{eq:3-9}
\end{align}
It is straightforward to check that for fixed $z>0$,
\begin{equation} \label{eq:3-10}
\sup_{0<\l x<z} x^{\l r_{1}}(z-\l x)^{(1-\l)r_{2}}
= {r_{1}^{\l r_{1}}\big((1-\l)r_{2}\big)^{(1-\l)r_{2}}
  z^{\l r_{1}+(1-\l)r_{2}} \over \big(\l r_{1}+(1-\l)r_{2}
  \big)^{\l r_{1}+(1-\l)r_{2}}}.
\end{equation}
Apply Equation (\ref{eq:3-10}) with $r_{1}=kt_{1}, r_{2}=
kt_{2}$ to Equation (\ref{eq:3-9}) to obtain
\begin{equation}
\ell_{u}(t_{1})^{\l} \ell_{u}(t_{2})^{1-\l} \geq
\left({(\l t_{1}+(1-\l)t_{2})^{\l t_{1}+(1-\l)t_{2}}
\over t_{1}^{\l t_{1}} t_{2}^{(1-\l)t_{2}}}\right)^{k}
\inf_{z>0} {u(z^{k}) \over z^{\l kt_{1}+(1-\l)kt_{2}}}.
      \notag
\end{equation}
But the last infimum is nothing but $\ell_{u}(\l t_{1}+
(1-\l)t_{2})$. Hence we have proved that
\begin{equation}
\ell_{u}(t_{1})^{\l} \ell_{u}(t_{2})^{1-\l} \geq
\left({(\l t_{1}+(1-\l)t_{2})^{\l t_{1}+(1-\l)t_{2}}
\over t_{1}^{\l t_{1}} t_{2}^{(1-\l)t_{2}}}\right)^{k}
\ell_{u}(\l t_{1}+(1-\l)t_{2}).  \notag
\end{equation}
This inequality shows that $\ell_{u}(t)t^{kt}$ is log-convex.

Conversely, suppose $\ell_{u}(t)t^{kt}$ is log-convex.
We can carry out similar calculations as above backward to
show that $u$ is (log, $x^{k}$)-convex.
\end{pf}

\begin{theorem} \label{lem:a-1}
Let $u\in C_{+, \log}$ be (log, $x^{k}$)-convex, $k>0$. Then
\begin{equation} \label{eq:a-1}
\ell_{u}(n) \ell_{u}(m) \leq \ell_{u}(0) 2^{k(n+m)}
  \ell_{u}(n+m), \qquad \forall n, m \geq 0.
\end{equation}
\end{theorem}

\noindent
{\it Remark.} Let $u\in C_{+, \log}$. By Theorem \ref{lem:3-2}
the sequence $\{\ell_{u}(n)\}_{n=0}^{\infty}$ is log-concave.
Then we can apply Theorem 2 (b) in \cite{akk1} with $\a(n)=
n!\ell_{u}(n)/\ell_{u}(0)$ to get
\begin{equation} \label{eq:a-2}
\ell_{u}(0) \ell_{u}(n+m) \leq \ell_{u}(n) \ell_{u}(m),
  \qquad \forall n, m \geq 0.
\end{equation}

\begin{pf}
By Lemma \ref{lem:3-5} the sequence $\{\ell_{u}(n) n^{kn}
\}_{n=0}^{\infty}$ is log-convex (here $0^{0}=1$ by
convention.) Then apply Theorem 2 (a) in \cite{akk1} with
$\a(n)= \ell_{u}(n) n^{kn}/\ell_{u}(0)$ to get
\begin{equation} \label{eq:a-3}
\ell_{u}(n)\ell_{u}(m) \leq \ell_{u}(0) \ell_{u}(n+m)
\left({(n+m)^{n+m} \over n^{n} m^{m}}\right)^{k}.
\end{equation}
Let $A=(n+m)^{n+m}/(n^{n}m^{m})$ and $x=n/(n+m)$. Then it
is easily checked that
\begin{equation}
{1\over n+m} \log A = -x \log x - (1-x)\log (1-x).
\end{equation}
But the maximum of the function $-x\log x-(1-x)\log (1-x)$
for $x\in (0, 1)$ obviously occurs at $x=1/2$ with a value
of $\log 2$. Hence
\begin{equation} \label{eq:a-4}
{(n+m)^{n+m} \over n^{n} m^{m}} = A \leq 2^{n+m}.
\end{equation}
Thus Equations (\ref{eq:a-3}) and (\ref{eq:a-4}) yield
Equation (\ref{eq:a-1}).
\end{pf}

\noindent
{\bf Inverse Legendre transform}

\smallskip
In view of Lemma \ref{lem:3-6}, we can define the inverse
Legendre transform as follows. Let $f$ be a positive
continuous function on $[0, \infty)$ such that
$\lim_{t\to\infty} f(t)^{1/t}=0$ or equivalently
$\lim_{t\to\infty} t^{-1} \log f(t) = -\infty$. Then we
define
\begin{equation} \label{eq:ilt}
\t_{f}(r) = \sup_{t\geq 0} f(t) r^{t}, \qquad r\geq 0.
\end{equation}
Suppose $u\in C_{+, \log}$ is (log, exp)-convex. By Theorems
\ref{lem:3-2} and \ref{lem:3-4}, $\t_{\ell_{u}}$ is defined.
Moreover, by Lemma \ref{lem:3-6}, we have
\begin{equation}
\t_{\ell_{u}}(r) = u(r), \qquad \forall r\geq 0. \notag
\end{equation}
Hence $\t_{\ell_{u}} = u$ for any (log, exp)-convex
funtion $u$ in $C_{+, \log}$.

On the other hand, let $f$ be a positive continuous
function on $[0, \infty)$ satisfying the conditions:
\begin{itemize}
\item[(1)] $\lim_{t\to\infty} f(t)^{1/t}=0$,
\item[(2)] $f$ is decreasing for large $t$,
\item[(3)] $f$ is log-concave.
\end{itemize}
We can carry out similar calculations as before to show
that $\ell_{\t_{f}} = f$. Therefore, $\t$ is the inverse
Legendre transform.

\medskip
Now, we come to questions (ii) and (iii) related to Equation
(\ref{eq:3-2}) in the beginning of this section. Note that
the coefficient of $x^{n}$ in Equation (\ref{eq:3-2}) is
$\ell_{u}(n)$. Hence the new function that we mentioned in
question (ii) is the series $\sum_{n=0}^{\infty} \ell_{u}(n)
x^{n}$. Since we will often refer to this function we give
it a name.

\begin{definition}
Let $u\in C_{+, \log}$ and $\lim_{n\to\infty}
\ell_{u}(n)^{1/n}=0$. The {\em $L$-function} of $u$ is
defined to be the function
\begin{equation} \label{eq:3-13}
\cl_{u}(r) = \sum_{n=0}^{\infty} \ell_{u}(n) r^{n},
  \qquad   r\geq 0.
\end{equation}
\end{definition}

Note that $\cl_{u}$ is an entire function. Let $u \in
C_{+, \log}$ be (log, exp)-convex. Then (i) by Theorem
\ref{lem:3-4} $\cl_{u}$ is defined and (ii) by Theorem
\ref{lem:3-2} the sequence $\{\ell_{u}(n)\}_{n=0}^{\infty}$
of coefficients in $\cl_{u}$ is log-concave.

\begin{lemma} \label{lem:a-2}
Let $u\in C_{+, \log}$ be (log, $x^{k}$)-convex, $k>0$.
Then
\begin{equation} \label{eq:a-5}
r\cl_{u}(r) \leq {\ell_{u}(0) \over \ell_{u}(1)}
  \cl_{u}(2^{k} r),  \qquad \forall r\geq 0.
\end{equation}
\end{lemma}

\begin{pf}
By Theorem \ref{lem:a-1} with $m=1$ we have
\begin{equation}
\ell_{u}(n)
\leq {\ell_{u}(0) \over \ell_{u}(1)} 2^{k(n+1)}
  \ell_{u}(n+1). \notag
\end{equation}
Hence for any $r\geq 0$,
\begin{equation}
r\cl_{u}(r)  = \sum_{n=0}^{\infty} \ell_{u}(n) r^{n+1}
\leq {\ell_{u}(0) \over \ell_{u}(1)} \sum_{n=0}^{\infty}
    2^{k(n+1)} \ell_{u}(n+1) r^{n+1}
\leq {\ell_{u}(0)\over\ell_{u}(1)} \cl_{u}(2^{k}r).\notag
\end{equation}
\end{pf}

\begin{theorem} \label{thm:3-1}
(1) Let $u\in C_{+, \log}$ be (log, exp)-convex. Then its
$L$-function $\cl_{u}$ is also (log, exp)-convex and for
any constant $a>1$,
\begin{equation} \label{eq:3-14}
  \cl_{u}(r) \leq {ea \over \log a}\,u(ar), \qquad
\forall r\geq 0.
\end{equation}

\noindent
(2) Let $u\in C_{+, \log}$ be increasing and
(log, $x^{k}$)-convex, $k>0$. Then there exists a constant
$C$, independent of $k$, such that
\begin{equation}  \label{eq:a-6}
u(r) \leq C \cl_{u}(2^{k}r), \qquad \forall r\geq 0.
\end{equation}
\end{theorem}

\noindent
{\it Remarks.} (a) From the proof below the constant
$C$ is given as follows. Note that if $u$ is increasing
and (log, $x^{k}$)-convex for some $k>0$, then by
Proposition \ref{pro:2-1} $u$ is (log, exp)-convex.
Hence by Theorem \ref{lem:3-4} its Legendre transform
$\ell_{u}(t)$ is decreasing for large $t$. Let $n_{0}$
be a natural number such that $\ell_{u}(t)$ is
decreasing for $t\geq n_{0}$. The constant $C$ is
given by
\begin{equation}
C = \max\left\{{u(1) \over \ell_{u}(0)}, {\ell_{u}(0)
\over \ell_{u}(1)}, {u(1) \over \ell_{u}(n_{0}+1)}
\right\}.  \notag
\end{equation}

\noindent
(b) If $u\in C_{+, \log}$ is increasing and
(log, $x^{k}$)-convex for some $k>0$, then we can combine
the inequalities in Equations (\ref{eq:3-14}) and
(\ref{eq:a-6}) together to get
\begin{equation} \label{eq:a-7}
{1\over C} u(2^{-k}r) \leq \cl_{u}(r) \leq {ea \over
\log a}\,u(ar), \qquad \forall r\geq 0.
\end{equation}

\begin{pf}
To prove the inequality in Equation (\ref{eq:3-14}), note
that from Equation (\ref{eq:3-8}) we have
\begin{equation}
\ell_{u}(s) \leq \r(t)^{t-s} \ell_{u}(t), \qquad
\forall s, t\geq 0,  \notag
\end{equation}
where $\r(t)$ is given in Lemma \ref{lem:3-3}. Hence for any
fixed $t\geq 0$,
\begin{align}
\cl_{u}(r) & = \sum_{n=0}^{\infty} \ell_{u}(n) r^{n} \notag\\
        & \leq \sum_{n=0}^{\infty} \r(t)^{t-n}
          \ell_{u}(t) r^{n} \notag  \\
     & = \r(t)^{t} \ell_{u}(t) \sum_{n=0}^{\infty}
          \big(r\r(t)^{-1}\big)^{n}.  \notag
\end{align}
Thus for $0<r<\r(t)$ we have
\begin{equation}
\cl_{u}(r) \leq \r(t)^{t} \ell_{u}(t) \big(1-r\r(t)^{-1}
   \big)^{-1}.   \notag
\end{equation}
Use this inequality to get
\begin{align}
\ell_{\cl_{u}}(t) & = \inf_{r>0} {\cl_{u}(r) \over r^{t}}
   \leq \inf_{0<r<\r(t)} {\cl_{u}(r) \over r^{t}} \notag\\
   & \leq \inf_{0<r<\r(t)} {\r(t)^{t} \ell_{u}(t)
     \big(1-r\r(t)^{-1}\big)^{-1} \over r^{t}} \notag \\
   & = \r(t)^{t} \ell_{u}(t) \inf_{0<r<\r(t)}
       \Big(r^{t}\big(1-r\r(t)^{-1}\big)\Big)^{-1} \notag
\end{align}
But it is easily checked that
\begin{equation}
\sup_{0<r<\r(t)} r^{t}\big(1-r\r(t)^{-1}\big) =
{t^{t} \r(t)^{t} \over (t+1)^{t+1}}. \notag
\end{equation}
Therefore, for any $t>0$, we have
\begin{equation}
\ell_{\cl_{u}}(t) \leq \ell_{u}(t) {(t+1)^{t+1} \over
  t^{t}}.  \notag
\end{equation}
Now, for any constant $a>1$, the inequality $(t+1)^{t+1}/
t^{t} \leq (ea/\log a) a^{t}$ holds for all $t>0$. Hence
\begin{equation}
\ell_{\cl_{u}}(t) \leq {ea \over \log a}\, a^{t}
     \ell_{u}(t), \qquad \forall t>0.  \notag
\end{equation}
By Lemma \ref{lem:3-1} $\ell_{\t_{a}u}(t) = a^{t}
\ell_{u}(t)$. Thus $\ell_{\cl_{u}}(t) \leq (ea/\log a)
\ell_{\t_{a}u}(t)$ for all $t\geq 0$. Then apply Lemma
\ref{lem:3-7}(2) to conclude that
\begin{equation}
\cl_{u}(r) \leq {ea \over \log a}\, u(ar),
     \qquad \forall r\geq 0.  \notag
\end{equation}

Now, we prove the inequality in Equation (\ref{eq:a-6}).
First suppose $0\leq r\leq 1$. Since $u$ and $\cl_{u}$ are
increasing functions and $\cl_{u}(r)\geq \ell_{u}(0)$,
we get
\begin{equation} \label{eq:a-8}
u(r) \leq u(1) \leq {u(1) \over \ell_{u}(0)}\cl_{u}(r)
\leq {u(1) \over \ell_{u}(0)}\cl_{u}(2^{k}r), \qquad
\forall r\in [0, 1].
\end{equation}

Before we consider $r>1$, let us note that by Proposition
\ref{pro:2-1} the function $u$, being increasing and
(log, $x^{k}$)-convex, is also (log, exp)-convex. Hence by
Theorem \ref{lem:3-4} $\ell_{u}(t)$ is decreasing for large
$t$. Let $n_{0}$ be a natural number such that
$\ell_{u}(t)$ is decreasing for $t\geq n_{0}$.

Let $r>1$ be fixed. By Lemma \ref{lem:3-6}, we have $u(r) =
\sup_{t\geq 0} \ell_{u}(t) r^{t}$. Hence there exists
$\tau =\tau(r)\geq 0$ such that
\begin{equation}
u(r) = \ell_{u}(\tau) r^{\tau}.  \notag
\end{equation}
Let $j=j(r)$ be the integer such that $j\leq \tau<j+1$.

\smallskip\noindent
{\sl Case 1:} $j\geq n_{0}$. In this case we have $u(r)\leq
\ell_{u}(j) r^{j+1}$ and so by Lemma \ref{lem:a-2}
\begin{equation} \label{eq:a-9}
u(r) \leq r\ell_{u}(j) r^{j} \leq r\cl_{u}(r) \leq
  {\ell_{u}(0) \over \ell_{u}(1)} \cl_{u}(2^{k} r).
\end{equation}

\smallskip\noindent
{\sl Case 2:} $j < n_{0}$. In this case, we use the fact
that $u(1)=\sup_{t\geq 0} \ell_{u}(t)$ to get
\begin{equation} \label{eq:a-10}
u(r) \leq u(1) r^{\tau} \leq u(1) r^{n_{0}+1} \leq
{u(1) \over \ell_{u}(n_{0}+1)} \cl_{u}(r) \leq
{u(1) \over \ell_{u}(n_{0}+1)} \cl_{u}(2^{k}r).
\end{equation}
Let $C=\max\{u(1)/\ell_{u}(0), \,\ell_{u}(0)/\ell_{u}(1),
\, u(1)/\ell_{u}(n_{0}+1)\}$. We can put Equations
(\ref{eq:a-8}), (\ref{eq:a-9}), and (\ref{eq:a-10})
together to get Equation (\ref{eq:a-6}).
\end{pf}

Now, observe that the inequalities in Equation
(\ref{eq:a-7}) are similar to those in Equation
(\ref{eq:k}). Thus the functions $u$ and $\cl_{u}$ are
what we called ``equivalent'' in the beginning of this
section. We now make this concept a formal definition.

\begin{definition} \label{def:eq}
Two positive functions $u$ and $v$ on $[0, \infty)$ are
called {\em equivalent} if there exist positive constants
$c_{1}, c_{2}, a_{1}, a_{2}$ such that
\begin{equation}
c_{1}u(a_{1}r) \leq v(r) \leq c_{2}u(a_{2}r), \qquad
\forall r\in [0, \infty). \notag
\end{equation}
\end{definition}

Suppose $u\in C_{+, \log}$ is increasing and
(log, $x^{2}$)-convex. Then by Theorem \ref{thm:3-1}
the function $u$ is equivalent to its $L$-function
$\cl_{u}$. Note that $\cl_{u}$ is (log, exp)-convex and
entire with positive coefficients. Moreover, Equation
(\ref{eq:a-6}) implies that $\cl_{u}\in C_{+, \log}$. Hence
we can state that each increasing (log, $x^{2}$)-convex
function in $C_{+, \log}$ is equivalent to a
(log, exp)-convex entire function with positive coefficients
in $C_{+, \log}$.

\begin{example}
Consider the function $u(r)=\exp\big[(1+\b)r^{1/(1+\b)}\big],
\,0\leq \b<1$. Obviously, $u \in C_{+, \log}$ is increasing
and (log, $x^{2}$)-convex. Its Legendre transform is easily
checked to be
\[ \ell_{u}(n)= \left\{
\begin{array}{ll}
\!\Big({e \over n}\Big)^{(1+\b)n}, & \mbox{if $n\geq 1$;} \\
1, & \mbox{if $n=0$}.
\end{array}\right.
\]
Hence the $L$-function of $u$ is given by
\begin{equation} \label{eq:3-15}
\cl_{u}(r) = \sum_{n=0}^{\infty} \Big({e\over n}
  \Big)^{(1+\b)n} r^{n},
\end{equation}
where $0^{0}=1$ by convention. We can use the Stirling
formula (see p.~357 in \cite{kuo96}) to get the inequalities
\begin{equation} \label{eq:3-16}
{1\over n!} \leq \Big({e\over n}\Big)^{n} \leq {e\,2^{n/2}
\over n!}, \qquad \forall n\geq 0.
\end{equation}
It follows from Equations (\ref{eq:3-15}) and (\ref{eq:3-16})
that
\begin{equation}
G_{1/\a}(r) \leq \cl_{u}(r) \leq e^{1+\b} G_{1/\a}\big(
2^{(1+\b)/2} r\big), \qquad r\geq 0, \notag
\end{equation}
where $G_{1/\a}(r)=\sum_{n=0}^{\infty} (n!)^{-(1+\b)} r^{n}$
as defined in Equation (\ref{eq:e}). Thus $\cl_{u}$ and
$G_{1/\a}$ are equivalent. On the other hand, by Theorem
\ref{thm:3-1}, $u$ is equivalent to $\cl_{u}$. Hence we
conclude that $u$ and $G_{1/\a}$ are equivalent.

On the other hand, consider the function $v(r)=\exp\big[
(1-\b)r^{1/(1-\b)}\big]$. By a similar argument as above
we can show that $v$ and the function $G_{\a}$ defined
in Equation (\ref{eq:e}) are equivalent. Note that
the functions $u$ and $v$ are nothing but $\wt G_{1/\a}$
and $\wt G_{\a}$, respectively, in Equation (\ref{eq:1-13}).
Thus the equivalence of $\wt G_{1/\a}$ and $\wt G_{\a}$ to
$G_{1/\a}$ and $G_{\a}$, respectively, has been proved
without using the inequalities in Equations (\ref{eq:c})
and (\ref{eq:d}) (cf.~(Q2) in Section \ref{sec:1}.)
\end{example}

At the end of this section we define the equivalence of
two sequences and state a simple fact which will be
convenient for future reference.

\begin{definition}
Two sequences $\{a(n)\}$ and $\{b(n)\}$ of nonnegative
numbers are said to be {\em equivalent} if there exist
positive constants $K_{1}, K_{2}, c_{1}, c_{2}$ such that
\begin{equation}  \label{eq:3-17}
K_{1}c_{1}^{n} a(n) \leq b(n) \leq K_{2} c_{2}^{n} a(n),
\qquad \forall n.
\end{equation}
\end{definition}

Let $f(r)$ and $g(r)$ be positive functions on $[0, \infty)$.
We want to point out that the equivalence of functions $f$
and $g$ (in the sense of Definition \ref{def:eq}) is quite
different from the equivalence of sequences $\{f(n)\}$ and
$\{g(n)\}$. Moreover, suppose $u(r)=\sum_{n=0}^{\infty} u_{n}
r^{n}$ is an entire function with $u_{n}>0$ and $\{u_{n}\}$
being log-concave. Then by Equation (\ref{eq:cks}) the
sequences $\{u_{n}\}$ and $\{\ell_{u}(n)\}$ are equivalent.

\begin{lemma} \label{lem:3-8}
Suppose $\{a(n)\}$ and $\{b(n)\}$ are equivalent sequences
of nonnegative numbers such that $a(n)^{1/n} \to 0$ or
$b(n)^{1/n} \to 0$ as $n\to\infty$. Then the functions
$A(r)=\sum_{n=0}^{\infty} a(n) r^{n}$ and
$B(r)=\sum_{n=0}^{\infty} b(n) r^{n}$ defined on
$[0, \infty)$ are equivalent.
\end{lemma}

\section{Dual Legendre function} \label{sec:4}

In this section we will develop a crucial machinery for the
next section and the application to white noise analysis in
the forthcoming paper \cite{akk3}.

We will think of the exponential generating function
$G_{1/\a}$ as $\cl_{u}$ for some
$u$. Equivalently, the sequence $\{\a(n)\}$ and the function
$u$ are related by the Legendre transform as follows:
\begin{equation} \label{eq:4-1}
\ell_{u}(n) = {1 \over n! \a(n)}.
\end{equation}
In that case the exponential function $G_{\a}$ is given by
\begin{equation} \label{eq:j}
G_{\a}(r)=\sum_{n=0}^{\infty} {\a(n) \over n!}\, r^{n}
= \sum_{n=0}^{\infty} {1 \over \ell_{u}(n) (n!)^{2}}\,
   r^{n}.
\end{equation}
But by Equation (\ref{eq:3-16}) the sequences $\{n!\}$ and
$\{(n/e)^{n}\}$ are equivalent. Hence by Lemma \ref{lem:3-8}
$G_{\a}$ is equivalent to the function defined by the series
\begin{equation} \label{eq:n}
\sum_{n=0}^{\infty} {e^{2n} \over \ell_{u}(n) n^{2n}}\,
r^{n}.
\end{equation}
A good way to understand this new function is to regard it
as $\cl_{v}$ for some $v$, i.e., we need to find $v$ such
that
\begin{equation} \label{eq:m}
\ell_{v}(t) = {e^{2t} \over \ell_{u}(t) t^{2t}}, \qquad
  t \geq 0,
\end{equation}
where $0^{0}=1$ by convention. The function $v$, defined
as $u^{*}$ in Definition \ref{def:4-1} below, belongs to
$C_{+, \log}$. Moreover, it is (log, exp)-convex by
Proposition \ref{pro:2-1} and Lemma \ref{lem:4-1} below.
Hence we can apply Lemma \ref{lem:3-6} to get
\begin{equation} \label{eq:4-a}
v(r) = \sup_{t\geq 0} \ell_{v}(t) r^{t} =
\sup_{t>0} {e^{2t} r^{t} \over \ell_{u}(t) t^{2t}}.
\end{equation}
Then use the definition of the Legendre transform to show
\begin{equation} \label{eq:4-2}
v(r)=\sup_{t, s>0} {e^{2t} r^{t} s^{t} \over u(s) t^{2t}}
= \sup_{s>0} {1\over u(s)} \sup_{t>0} {(e^{2}rs)^{t}
  \over t^{2t}}.
\end{equation}
But it can be easily checked that for $a>0$,
\begin{equation} \label{eq:4-3}
\sup_{t>0} {a^{t} \over t^{2t}} = e^{2\sqrt{a}/e}.
\end{equation}
Put Equation (\ref{eq:4-3}) with $a=e^{2}rs$ into Equation
(\ref{eq:4-2}) to conclude that
\begin{equation} \label{eq:4-4}
v(r) = \sup_{s>0} {e^{2\sqrt{rs}} \over u(s)}.
\end{equation}
This equation suggests a new transform and raises a question
of finding $u$ for which this new transform can be defined.

\smallskip
\noindent
{\bf Notation.} Let $C_{+, j}, j>0,$ denote the set of all
positive continuous functions $u$ on $[0, \infty)$ satisfying
the condition
\begin{equation} \label{eq:4-5}
\lim_{r\to\infty} {\log u(r) \over r^{j}} = \infty.
\end{equation}

\smallskip
We will mostly be concerned with the set $C_{+, 1/2}$ because
the right-hand side of Equation (\ref{eq:4-4}) exists for
all $r\geq 0$ when $u\in C_{+, 1/2}$. On the other hand,
observe that $C_{+, j}\subset C_{+, \log}$ for all $j>0$.

\begin{definition} \label{def:4-1}
The {\em dual Legendre function} $u^{*}$ of $u \in
C_{+, 1/2}$ is defined to be the function
\begin{equation} \label{eq:4-6}
u^{*}(r) = \sup_{s>0} {e^{2\sqrt{rs}} \over u(s)},
\qquad r\geq 0.
\end{equation}
\end{definition}

\noindent
{\it Remark.}
In \cite{ghor} by Gannoun et al., 
they adopted the relation:
$$
	\t(r)^{*}:=\sup_{s\geq 0}\{sr-\t(s)\}, \quad r\geq 0.
$$
Hence 
\begin{equation}
	\log u(s)=\t(\sqrt{2s}),\quad
	\log u(r)^{*}=\t(\sqrt{2r})^{*}
\end{equation}
hold.

\begin{example}
For the function $u(r)=e^{r}$, we have $u^{*}(r)=e^{r}$.
This is the case for the Hida-Kubo-Takenaka space.
\end{example}

\begin{example}
For the function $u(r)=\exp\big[(1+\b)r^{1/(1+\b)}\big]$,
we can easily check that $u^{*}(r)=\exp\big[(1-\b)
r^{1/(1-\b)}\big]$. This is the case for the
Kondratiev-Streit space.
\end{example}

\begin{example}
Let $u(r)=\exp[e^{x}]$. To find $u^{*}(r)$ we need to find
the maximum of the function $2\sqrt{rs} -e^{s}$. The critical
point $s_{0}$ of this function satisfies the equation
\begin{equation}
\sqrt{r} = \sqrt{s}\,e^{s}. \notag
\end{equation}
Obviously, we have $\lim_{r\to\infty} s_{0}=\infty$. Hence
$s_{0}\sim \log\sqrt{r}\,$ for large $r$ and so
\begin{equation}
\sup_{s>0} \big(2\sqrt{rs} - e^{s}\big) = 2\sqrt{rs_{0}} -
e^{s_{0}} = 2\sqrt{rs_{0}} - {\sqrt{r} \over \sqrt{s_{0}}}
\sim 2\sqrt{rs_{0}} \sim 2\sqrt{r\log\sqrt{r}}. \notag
\end{equation}
Thus although we cannot find the exact form of $u^{*}$,
the function $u^{*}$ is equivalent to the function
$\exp\big[2\sqrt{r\log\sqrt{r}}\,\big]$. In general let
\begin{equation}
u(r)=\exp_{k}(r)=\exp(\exp(\cdots(\exp(r)))), \qquad
\text{$k$-th iteration}.  \notag
\end{equation}
Its dual Legendre function $u^{*}$ is equivalent to the
function
\begin{equation}
\exp\left[2\,\sqrt{r\log_{k-1}\sqrt{r}}\right], \notag
\end{equation}
where $\log_{j}$ is defined by
\begin{equation}
\log_{1}(r)=\log(\max\{r, e\}), \quad \log_{j}(r) =
  \log_{1}(\log_{j-1}(r)), \quad j\geq 2. \notag
\end{equation}
This example is for the Gel'fand triple associated with the
Bell numbers in the paper by Cochran et al.~\cite{cks}.
\end{example}

\begin{lemma} \label{lem:4-1}
Let $u\in C_{+, 1/2}$. Then its dual Legendre function
$u^{*}$ belongs to $C_{+, 1/2}$ and is an increasing
(log, $x^{2}$)-convex function.
\end{lemma}

\begin{pf}
From the definition of $u^{*}(r)$ we have $\log u^{*}(r)
\geq 2\sqrt{rs} - \log u(s)$ for any $s>0$. Hence
\begin{equation}
{\log u^{*}(r) \over \sqrt{r}} \geq {2\sqrt{rs} - \log u(s)
\over\sqrt{r}}=2\sqrt{s}-{\log u(s) \over \sqrt{r}}. \notag
\end{equation}
This implies that
\begin{equation}
\liminf_{r\to\infty} {\log u^{*}(r) \over \sqrt{r}}
\geq 2\sqrt{s}, \qquad \forall s>0.  \notag
\end{equation}
Therefore, $\lim_{r\to\infty} \log u^{*}(r)/\sqrt{r} =
\infty$, which shows that $u^{*} \in C_{+, 1/2}$. To show
that $u^{*}(r)$ is increasing, let $r_{1}<r_{2}$. Note that
there exists some $s_{1}>0$ such that $u^{*}(r_{1})=
e^{2\sqrt{r_{1}s_{1}}}/u(s_{1})$. Hence
\begin{equation}
u^{*}(r_{1}) = {e^{2\sqrt{r_{1}s_{1}}} \over u(s_{1})}
  \leq {e^{2\sqrt{r_{2}s_{1}}} \over u(s_{1})}
  \leq u^{*}(r_{2}).  \notag
\end{equation}
To show that $u^{*}(r)$ is (log, $x^{2}$)-convex, let
$r_{1}, r_{2}\geq 0$, and $0\leq \l\leq 1$. Then
\begin{align}
u^{*}\big((\l r_{1}+(1-\l)r_{2})^{2}\big)
  & = \sup_{s>0} {e^{2(\l r_{1}+(1-\l)r_{2})\sqrt{s}}
      \over u(s)}   \notag  \\
  & \leq \left(\sup_{s_{1}>0} {e^{2r_{1}\sqrt{s_{1}}} \over
    u(s_{1})}\right)^{\l} \left(\sup_{s_{2}>0}
    {e^{2r_{2}\sqrt{s_{2}}} \over u(s_{2})}
     \right)^{1-\l}  \notag  \\
  & = u^{*}(r_{1}^{2})^{\l}\,
        u^{*}(r_{2}^{2})^{1-\l}. \notag
\end{align}
Thus $u^{*}(r)$ is (log, $x^{2}$)-convex.
\end{pf}

\begin{theorem} \label{thm:4-2}
Let $u\in C_{+, 1/2}$ be (log, $x^{2}$)-convex. Then
the Legendre transform of $u^{*}$ is given by
\begin{equation} \label{eq:4-7}
\ell_{u^{*}}(t) = {e^{2t} \over \ell_{u}(t) t^{2t}}.
\end{equation}
\end{theorem}

\begin{pf}
Note that $u^{*} \in C_{+, \log}$ since $u^{*}\in
C_{+, 1/2}$ by Lemma \ref{lem:4-1} and $C_{+, j}\subset
C_{+, \log}$ for all $j>0$. Hence the Legendre transform
$\ell_{u^{*}}$ is defined. By assumption $u$ is
(log, $x^{2}$)-convex and so by Lemma \ref{lem:3-5} the
function $\ell_{u}(t) t^{2t}$ is log-convex. Hence
$\big(\ell_{u}(t) t^{2t}\big)^{-1}$ is log-concave. Since
$e^{2t}$ is also log-concave, we see that the function
\begin{equation}
w(t) = {e^{2t} \over \ell_{u}(t) t^{2t}} \notag
\end{equation}
is log-concave. Note that $\ell_{u}(t)^{1/t} t^{2}$
increases to $\infty$ as $t\to\infty$ since the function
$u$ is (log, $x^{2}$)-convex. Hence the inverse Legendre
transform $\t$ in Equation (\ref{eq:ilt}) is defined at
$w$ by
\begin{equation} \label{eq:4-10}
\t_{w}(r) = \sup_{t\geq 0} {e^{2t} r^{t} \over
\ell_{u}(t) t^{2t}}.
\end{equation}
Moreover, $\ell_{\t_{w}} = w$. On the other hand, from the
motivation for the dual Legendre function in Equations
(\ref{eq:4-a}) (\ref{eq:4-4}) (\ref{eq:4-6}) we have
\begin{equation} \label{eq:4-11}
u^{*}(r) = \sup_{t\geq 0} {e^{2t} r^{t} \over
\ell_{u}(t) t^{2t}}.
\end{equation}
It follows from Equations (\ref{eq:4-10}) and (\ref{eq:4-11})
that $\t_{w} = u^{*}$. But we also have $\ell_{\t_{w}} = w$.
Hence $\ell_{u^{*}} = w$ and the theorem is proved.
\end{pf}

\noindent
{\it Remark.}
Let $u\in C_{+, 1/2}$ be (log, $x^{2}$)-convex. Suppose
$u$ is increasing on the interval $[r_{0}, \infty)$. Then
$(u^{*})^{*}(r) = u(r)$ for all $r\geq r_{0}$.
Observe that if $u$ is an increasing (log, $x^{2}$)-convex
function in $C_{+, 1/2}$, then we have $(u^{*})^{*} = u$.
Since we will not use this involution property elsewhere
in this paper, we skip the proof.

As we mentioned in the beginning of this section the
exponential generating function $G_{1/\a}$ 
is thought of as the $L$-function $\cl_{u}$ for
some function $u$. Then the corresponding exponential
generating function $G_{\a}$, expressed in terms of
$\ell_{u}(n)$'s, is given by the second series in Equation
(\ref{eq:j}). We give this series a name for future
reference.

\begin{definition}
Let $u\in C_{+, 1/2}$ and suppose $\lim_{n\to\infty}
\big(\ell_{u}(n)(n!)^{2}\big)^{-1/n}=0$. The
{\em $L^{\#}$-function} of $u$ is defined to be the function
\begin{equation}
\cl_{u}^{\#}(r) = \sum_{n=0}^{\infty} {1 \over \ell_{u}(n)
  (n!)^{2}}\, r^{n}, \qquad r\geq 0.  \notag
\end{equation}
\end{definition}

Note that $\cl_{u}^{\#}(r)$ is an entire function. It
follows from Theorem \ref{thm:4-2} and Equations
(\ref{eq:j}) and (\ref{eq:n}) that $\cl_{u}^{\#}(r)$ is
defined for any (log, $x^{2}$)-convex function $u$ in
$C_{+, 1/2}$.

\begin{theorem} \label{thm:4-3}
Let $u\in C_{+, 1/2}$ be (log, $x^{2}$)-convex. Then the
functions $\cl_{u^{*}}$ and $\cl_{u}^{\#}$ are equivalent.
\end{theorem}

\noindent
{\it Remark.} Let $u\in C_{+, 1/2}$. Then by Lemma
\ref{lem:4-1} its dual Legendre transform $u^{*}$ belongs
to $C_{+, 1/2}$ and is increasing and (log, $x^{2}$)-convex.
Hence we can apply Theorem \ref{thm:3-1} to $u^{*}$ to
conclude that the functions $u^{*}$ and $\cl_{u^{*}}$ are
equivalent. Therefore, under the assumption of the above
theorem, the functions $u^{*},\,\cl_{u^{*}},\,\cl_{u}^{\#}$
are all equivalent.

\begin{pf}
Note that the function $\cl_{u}^{\#}$ is the second series
in Equation (\ref{eq:j}). But from the discussion for
Equations (\ref{eq:n}) (\ref{eq:m}) (\ref{eq:4-a}) with $v=
u^{*}$, we see easily that $\cl_{u}^{\#}$ is equivalent to
the function $\sum_{n=0}^{\infty} \ell_{u^{*}}(n) r^{n}$,
which is exactly the function $\cl_{u^{*}}$.
\end{pf}

Below we list some facts concerning equivalent functions
and sequences. These facts can be easily checked by using
the previous results or the techniques in the proofs.

\begin{itemize}
\item[1.] If $u\in C_{+, \log}$ and $v$ is equivalent to $u$,
then $v\in C_{+, \log}$ and the sequences $\{\ell_{u}(n)\}$
and $\{\ell_{v}(n)\}$ are equivalent.
\item[2.] If $u, v \in C_{+, \log}$, $u$ is increasing and
(log, $x^{2}$)-convex, and the sequences $\{\ell_{u}(n)\}$
and $\{\ell_{v}(n)\}$ are equivalent, then the functions
$u$ and $v$ are equivalent.
\item[3.] If $u\in C_{+, \log}$, $u$ is increasing and
(log, $x^{2}$)-convex, and $u$ and $v$ are equivalent, then
the $L$-functions $\cl_{u}$ and $\cl_{v}$ are equivalent.
\item[4.] If $u\in C_{+, 1/2}$ and $v$ is equivalent to $u$,
then $v\in C_{+, 1/2}$ and the functions $u^{*}$ and $v^{*}$
are equivalent.
\item[5.] If $u\in C_{+, 1/2}$, $u$ is (log, $x^{2}$)-convex,
and $u$ and $v$ are equivalent, then the functions
$\cl_{u}^{\#}$ and $\cl_{v}^{\#}$ are equivalent.
\end{itemize}

Many properties of a function or sequence remain true for
equivalent functions or sequences. For convenience, we make
the following definition.

\begin{definition}
Let $P$ be a property of functions or sequences. A function
$u$ is said to be {\em nearly $P$} if there exists a $P$
function which is equivalent to $u$. A sequence $\{a(n)\}$
is said to be {\em nearly $P$} if there exists a $P$
sequence which is equivalent to $\{a(n)\}$.
\end{definition}

For example, a positive function $u$ is nearly
(log, exp)-convex if there exists a (log, exp)-convex
function which is equivalent to $u$. A positive sequence
$\{a(n)\}$ is nearly log-concave if there exists a
log-concave sequence which is equivalent to $\{a(n)\}$.

Here we list some results concerning functions and sequences
that are ``nearly'' something.

\begin{itemize}
\item[6.] Let $u, v\in C_{+, \log}$ be increasing and nearly
(log, $x^{2}$)-convex. Then the functions $u$ and $v$ are
equivalent if and only if the sequences $\{\ell_{u}(n)\}$
and $\{\ell_{v}(n)\}$ are equivalent.
\item[7.] Let $u, v\in C_{+, 1/2}$ be nearly
(log, $x^{2}$)-convex. Then $u$ and $v$ are equivalent if
and only if $u^{*}$ and $v^{*}$ are equivalent.
\item[8.] Let $u(r)=\sum_{n=0}^{\infty} u_{n}r^{n}$ and
$v(r)=\sum_{n=0}^{\infty} v_{n}r^{n}$ be entire functions
with $u_{n}, v_{n}>0$. Suppose $\{u_{n}\}$ and $\{v_{n}\}$
are nearly log-concave sequences. Then $\{u_{n}\}$ and
$\{v_{n}\}$ are equivalent if and only if $u$ and $v$ are
equivalent if and only if $\{\ell_{u}(n)\}$ and
$\{\ell_{v}(n)\}$ are equivalent.
\item[9.] If $u\in C_{+, \log}$, then the sequence
$\{\ell_{u}(n)\}$ is log-concave. On the other hand, if
$u\in C_{+, 1/2}$ is (log, $x^{2}$)-convex, then the
sequence $\big\{\big(\ell_{u}(n)(n!)^{2}\big)^{-1}\big\}$
is nearly log-concave.
\end{itemize}

We make two remarks about Item 9: (1) Let $\{b_{k}(n)\}$
be the Bell numbers of order $k$. It has been shown in
\cite{akk1} that $\{b_{k}(n)/n!\}$ is log-concave and
$\{b_{k}(n)\}$ is log-convex. Note that $\{b_{k}(n)\}$
being log-convex implies that $\{(b_{k}(n)n!)^{-1}\}$ is
log-concave.

\noindent
(2) The near log-concavity of the sequence
$\big\{\big(\ell_{u}(n)(n!)^{2}\big)^{-1}\big\}$ has been
shown in \cite{kubo} to be a necessary condition for the
characterization theorem of generalized functions in the
Gel'fand triple introduced by Cochran et al.~\cite{cks}.

\section{Growth order of holomorphic functions} \label{sec:5}

Recall that the $S$-transform $F$ of a generalized function
is a function on the complexification $\ce_{c}$ of $\ce$.
It is a holomorphic function on $\ce_{c}$ in the sense that
for any $\x, \y \in \ce_{c}$, the function $F(z\x+\y)$ is
an entire function of $z\in\spc$. Moreover, it satisfies
the growth conditions in Equations (\ref{eq:1-10}) and
(\ref{eq:1-11}) for generalized and test functions,
respectively. In this section we will study the
representation of holomorphic functions $F$ on $\ce_{c}$
satisfying the growth conditions in Equations (\ref{eq:1-10})
and (\ref{eq:1-11}) with $G_{\a}$ and $G_{1/\a}$ being
replaced by certain functions. The characterization theorems
will be given in our forthcoming papers \cite{akk3}.

\begin{lemma} \label{lem:5-1}
Let $u\in C_{+, \log}$. Suppose $F$ is a holomorphic function
on $\ce_{c}$ and there exist constants $K, a, p \geq 0$ such
that
\begin{equation} \label{eq:5-1}
|F(\x)| \leq K\,u\big(a|\x|_{-p}^{2}\big)^{1/2}, \qquad
\forall \x\in \ce_{c}.
\end{equation}
Let $q\in [0, p]$ be an integer such that $i_{p, q}$ is a
Hilbert-Schimidt operator. Then there exist functions
$f_{n}\in \ce_{q, c}^{\wh\otimes n}$ such that $F(\x)=
\sum_{n=0}^{\infty} \la f_{n}, \x^{\wh\otimes n} \ra$ and
\begin{equation} \label{eq:5-2}
|f_{n}|_{q}^{2} \leq K^{2}\big(a e^{2}\|i_{p, q}\|_{HS}^{2}
\big)^{n}\,\ell_{u}(n),
\end{equation}
where $\ell_{u}$ is the Legendre transform of $u$.
\end{lemma}

\begin{pf}
We follow the same argument as in the proof of Theorem 8.9
in \cite{kuo96}. Since $F$ is a holomorphic function on
$\ce_{c}$, it has the expansion
\begin{equation}
F(\x) = \sum_{n=0}^{\infty}J_{n}(\x, \x, \ldots, \x), \notag
\end{equation}
where $J_{n}$ is a symmetric $n$-linear functional on
$\ce_c\times\cdots \times\ce_c$ given by
\begin{equation}
J_n(\x_1, \ldots, \x_n) = {1\over n!}
{\partial\over \partial z_1}\cdots {\partial\over
\partial z_n} F(z_1\x_1+\cdots+z_n\x_n)
\Big|_{z_1=\cdots=z_n=0}. \notag
\end{equation}
Apply the Cauchy formula to show that
\begin{equation}
J_n(\x_1, \ldots, \x_n) = {1\over n!} {1\over (2\p i)^{n}}
\int_{|z_{1}|=r_{1}}\cdots \int_{|z_{n}|=r_{n}}
{F(z_{1}\x_{1}+\cdots+z_{n}\x_{n}) \over z_{1}^{2}\cdots
z_{n}^{2}}\,dz_{1}\cdots dz_{n}. \notag
\end{equation}
Let $R>0$. For nonzero $\x_{j}$'s, take $r_{j} =
R/|\x_{j}|_{-p}, \,\, 1\leq j\leq n$ and use the maximum
modulus principle to derive that
\begin{equation}
|J_n(\x_1, \ldots, \x_n)| \leq {1\over n!} {1\over R^{n}}
\left(\sup_{|\x|_{-p}=nR} |F(\x)|\right) |\x_{1}|_{-p}
\cdots |\x_{n}|_{-p}.  \notag
\end{equation}
Use the growth condition in Equation (\ref{eq:5-1}) to get
\begin{equation}
|J_n(\x_1, \ldots, \x_n)| \leq K {1\over n!} {u\big(an^{2}
R^{2}\big)^{1/2} \over R^{n}}\,|\x_{1}|_{-p}
\cdots |\x_{n}|_{-p}.  \notag
\end{equation}
This inequality holds for any $R>0$ and $\x_{1}, \ldots,
\x_{n} \in \ce_{p}'$. Let $an^{2}R^{2}=r$. Then
\begin{equation}
|J_n(\x_1, \ldots, \x_n)| \leq K {a^{n/2} n^{n} \over n!}
\left({u(r) \over r^{n}}\right)^{1/2} |\x_{1}|_{-p}
\cdots |\x_{n}|_{-p}.  \notag
\end{equation}
Now, take the infimum over $r>0$ to obtain
\begin{equation}
|J_n(\x_1, \ldots, \x_n)| \leq K {a^{n/2} n^{n} \over n!}
\,\ell_{u}(n)^{1/2}\,|\x_{1}|_{-p} \cdots |\x_{n}|_{-p}.
  \notag
\end{equation}
Then use the same argument as in the proof of Theorem 8.9 in
\cite{kuo96} to conclude that $J_n(\x_1, \ldots, \x_n) =
\la f_{n}, \x_{1}\otimes \cdots \x_{n}\ra$ with $f_{n} \in
\ce_{q, c}^{\wh\otimes n}, \, q\in [0, p]$, and
\begin{equation}
|f_{n}|_{q}^{2} \leq K^{2} a^{n}\left({n^{n} \over n!}
\right)^{2} \,\ell_{u}(n)\,\|i_{p, q}\|_{HS}^{2n}. \notag
\end{equation}
This inequality implies the one in Equation (\ref{eq:5-2})
because $n^{n} \leq n!e^{n}$.
\end{pf}

\begin{lemma} \label{lem:5-2}
Let $u\in C_{+, \log}$ be (log, exp)-convex. Suppose $F(\x)
=\sum_{n=0}^{\infty} \la f_{n}, \x^{\otimes n}\ra$ is a
holomorphic function on $\ce_{c}$ and there exist $K, a,
p\geq 0$ such that
\begin{equation}
|f_{n}|_{p} \leq K a^{n} \ell_{u}(n)^{1/2}, \qquad
\forall n\geq 0.  \notag
\end{equation}
Then for any $\x\in\ce_{c}$,
\begin{equation} \label{eq:5-3}
|F(\x)| \leq \sqrt{2}\,eK u\big(2ea^{2}|\x|_{-p}^{2}
\big)^{1/2}.
\end{equation}
\end{lemma}

\begin{pf}
By assumption we have $|\la f_{n}, \x^{\otimes n}\ra|
\leq Ka^{n} \ell_{u}(n)^{1/2} |\x|_{-p}^{n}$. Hence
\begin{align}
|F(\x)| & \leq \sum_{n=0}^{\infty} K a^{n}\ell_{u}(n)^{1/2}
  |\x|_{-p}^{n}  \notag  \\
  & = K\sum_{n=0}^{\infty}\left({1\over \sqrt{2}}\right)^{n}
     \left(\ell_{u}(n)^{1/2}\big(\sqrt{2}\,a |\x|_{-p}
     \big)^{n}\right)  \notag  \\
  & \leq K \sqrt{2} \left(\sum_{n=0}^{\infty} \ell_{u}(n)
    \big(2a^{2} |\x|_{-p}^{2}\big)^{n}\right)^{1/2} \notag\\
  & =  K \sqrt{2} \,\cl_{u}\big(2a^{2}
        |\x|_{-p}^{2}\big)^{1/2},  \label{eq:5-4}
\end{align}
where $\cl_{u}$ is the $L$-function of $u$. But by Theorem
\ref{thm:3-1} (1) with $a=e$ we have
\begin{equation} \label{eq:5-5}
\cl_{u}(r)\leq e^{2} u(er), \qquad r\geq 0.
\end{equation}
The conclusion in Equation (\ref{eq:5-3}) follows from
Equations (\ref{eq:5-4}) and (\ref{eq:5-5}).
\end{pf}

Now, let $u \in C_{+, \log}$ be a fixed function. Suppose
$F$ is an entire function on $\ce_{c}$ with the expansion
$F(\x)=\sum_{n=0}^{\infty} \la f_{n}, \x^{\otimes n}\ra$.
Being motivated by the norm given in \cite{cks} (with
$\ell_{u}(n)$ replacing $(n!\a(n))^{-1}$ as noted before,)
we define for each $p\geq 0$,
\begin{equation}
\|F\|_{u, p}=\left(\sum_{n=0}^{\infty} {1\over \ell_{u}(n)}
  \,|f_{n}|_{p}^{2}\right)^{1/2}. \notag
\end{equation}
Let $\ck_{u, p} = \{F; \|F\|_{u, p}< \infty\}$. Then
$\ck_{u, p}$ is a Hilbert space with norm $\|\cdot\|_{u, p}$.

On the other hand, being motivated by the work of Lee
\cite{lee} and Section 15.2 in the book \cite{kuo96},
we define $|\!|\!|F|\!|\!|_{u, p}$ for a holomorphic function $F$
on $\ce_{c}$ and for each $p\geq 0$ by
\begin{equation}
|\!|\!|F|\!|\!|_{u, p} = \sup_{\x\in\ce_{c}} |F(\x)|
  u\big(|\x|_{-p}^{2}\big)^{-1/2}.  \notag
\end{equation}
Let $\cg_{u, p}=\{F; |\!|\!|F|\!|\!|_{u, p}<\infty\}$. Then
$\cg_{u, p}$ is a Banach space with norm
$|\!|\!|\cdot|\!|\!|_{u, p}$.

\begin{theorem} \label{thm:5-1}
Let $u\in C_{+, \log}$. Suppose $p>q$ is such that the
inclusion mapping $i_{p, q}: \ce_{p} \to \ce_{q}$ is a
Hilbert-Schmidt operator with $\|i_{p, q}\|_{HS} \leq
e^{-1}$. Then
\begin{equation} \label{eq:5-6}
\|F\|_{u, q} \leq \left(1-e^{2}\|i_{p, q}\|_{HS}^{2}
  \right)^{-1/2} |\!|\!|F|\!|\!|_{u, p}, \qquad \forall
  F\in \cg_{u, p}.
\end{equation}
\end{theorem}

\noindent
{\bf Remark.} Conditions (a) and (b) stated in the beginning
of Section \ref{sec:1} imply that $\lim_{p\to\infty}
\|i_{p, q}\|_{HS}=0$ for any $q\geq 0$. Hence for any given
$q\geq 0$, there exists some $p>q$ such that
$\|i_{p, q}\|_{HS} \leq e^{-1}$. Therefore, it follows from
the theorem that for any $q\geq 0$, there exists $p>q$ such
that $\cg_{u, p} \subset \ck_{u, q}$ and the inclusion
mapping is continuous by Equation (\ref{eq:5-6}).

\begin{pf}
Suppose $F\in \cg_{u, p}$. Then we have
\begin{equation}
|F(\x)| \leq |\!|\!|F|\!|\!|_{u, p}\,u\big(|\x|_{-p}^{2}
  \big)^{1/2}, \qquad \forall \x\in\ce_{c}.
\end{equation}
Hence for $q$ as specified in the theorem, we can apply
Lemma \ref{lem:5-1} to show that $F(\x)=\sum_{n=0}^{\infty}
\la f_{n}, \x^{\otimes n}\ra$ with $f_{n}\in
\ce_{q, c}^{\wh\otimes n}$ and
\begin{equation}
|f_{n}|_{q}^{2} \leq |\!|\!|F|\!|\!|_{u, p}^{2}\left(e^{2}
\|i_{p, q}\|_{HS}^{2}\right)^{n} \ell_{u}(n). \notag
\end{equation}
Therefore,
\begin{align}
\|F\|_{u, q}^{2}
  & = \sum_{n=0}^{\infty} {1\over \ell_{u}(n)}
      |f_{n}|_{q}^{2}  \notag  \\
& \leq \sum_{n=0}^{\infty} {1\over \ell_{u}(n)}
  |\!|\!|F|\!|\!|_{u, p}^{2}\left(e^{2}\|i_{p, q}\|_{HS}^{2}
  \right)^{n} \ell_{u}(n)  \notag  \\
  & = \left(1-e^{2}\|i_{p, q}\|_{HS}^{2}
       \right)^{-1} |\!|\!|F|\!|\!|_{u, p}^{2}.  \notag
\end{align}
This proves the inequality in Equation (\ref{eq:5-6}).
\end{pf}

\begin{theorem} \label{thm:5-2}
Let $u\in C_{+, \log}$ be (log, exp)-convex. Then for any
$p\geq 1$, we have
\begin{equation} \label{eq:5-7}
|\!|\!|F|\!|\!|_{u, p-1} \leq \sqrt{e}\big(2\r^{2} \log 1/\r
\big)^{-1/2}\,\|F\|_{u, p}, \qquad \forall F\in \ck_{u, p},
\end{equation}
where the constant $\r$ is given in Condition (a) in the
beginning of Section \ref{sec:1}.
\end{theorem}

\noindent
{\bf Remark.}
It follows from Equation (\ref{eq:5-7}) that for any
$p\geq 1, \, \ck_{u, p} \subset \cg_{u, p-1}$ and the
inclusion mapping is continuous.

\begin{pf}
Let $F\in \ck_{u, p}$ and $p\geq 1$. Since $F(\x)=
\sum_{n=0}^{\infty} \la f_{n}, \x^{\wh\otimes n}\ra$,
we can derive that
\begin{align}
|F(\x)|
& \leq \sum_{n=0}^{\infty} |f_{n}|_{p} |\x|_{-p}^{n}
        \notag  \\
& = \sum_{n=0}^{\infty} \left({1\over \sqrt{\ell_{u}(n)}}
   \,|f_{n}|_{p}\right)\left(\sqrt{\ell_{u}(n)}\,
   |\x|_{-p}^{n}\right)  \notag \\
& \leq \left(\sum_{n=0}^{\infty} {1\over \ell_{u}(n)}
   \,|f_{n}|_{p}^{2}\right)^{1/2} \left(\sum_{n=0}^{\infty}
   \ell_{u}(n)\,|\x|_{-p}^{2n} \right)^{1/2} \notag \\
& = \|F\|_{u, p} \,\cl_{u}\big(|\x|_{-p}^{2}
     \big)^{1/2}.  \label{eq:5-8}
\end{align}
Note that $|\x|_{-p}\leq \r |\x|_{-p+1}$ and then apply
Theorem \ref{thm:3-1} (1) with $a=1/\r^{2}$ to get
\begin{equation} \label{eq:5-9}
\cl_{u}\big(|\x|_{-p}^{2}\big) \leq \cl_{u}\big(\r^{2}
|\x|_{-p+1}^{2}\big) \leq e \big(2\r^{2} \log 1/\r\big)^{-1}
\, u\big(|\x|_{-p+1}^{2}\big).
\end{equation}
Equations (\ref{eq:5-8}) and (\ref{eq:5-9}) imply the
inequality in Equation (\ref{eq:5-7}).
\end{pf}

Take a (log, exp)-convex function $u\in C_{+, \log}$.
Let $\ck_{u}$ and $\cg_{u}$ be the projective limits of
the families $\{\ck_{u, p}; p\geq 0\}$ and $\{\cg_{u, p};
p\geq 0\}$, respectively. By the remarks following each of
Theorems \ref{thm:5-1} and \ref{thm:5-2} we see that
$\ck_{u}=\cg_{u}$ and their respective topologies given by
$\{\|\cdot\|_{u, p}; p\geq 0\}$ and
$\{|\!|\!|\cdot|\!|\!|_{u, p}; p\geq 0\}$ coincide. In the
forthcoming paper we will study the corresponding spaces of
test and generalized functions and the characterization
theorems. \\
\bigskip
\begin{center}
{\appendixname}
\end{center}

Let $\{\a(n)\}_{n=0}^{\infty}$ be a sequence of positive
numbers. Let us extract the following list from \cite{akk4}:

\medskip
\begin{itemize}
\item[(A1)] $\a(0)=1$ and $\inf_{n\geq 0} \a(n) \s^{n}
>0$ for some $\s\geq 1$.
\medskip
\item[(A2)] $\lim_{n\to\infty} \left({\a(n) \over n!}
\right)^{1/n} =0$.
\medskip
\item[($\wt{\text{A}}2$)] $\lim_{n\to\infty} \left({1\over
  n!\a(n)}\right)^{1/n} =0$.
\smallskip
\item[(B1)] $\limsup_{n\to\infty} \left({n! \over \a(n)}
\inf_{r>0} {G_{\a}(r) \over r^{n}}\right)^{1/n} < \infty$.
\smallskip
\item[($\wt{\text{B}}1$)] $\limsup_{n\to\infty} \left(
n!\a(n) \inf_{r>0} {G_{1/\a}(r) \over r^{n}}\right)^{1/n}
  < \infty$.
\medskip
\item[(B2)] The sequence $\g(n)={\a(n) \over n!}, n\geq 0,$
is log-concave, i.e., for all $n\geq 0$,
\begin{equation}
\g(n) \g(n+2) \leq \g(n+1)^{2}. \notag
\end{equation}
\item[($\wt{\text{B}}2$)] The sequence $\left\{{1 \over
n!\a(n)}\right\}$ is log-concave.
\bigskip
\item[(B3)] The sequence $\{\a(n)\}$ is log-convex, i.e.,
for all $n\geq 0$,
\begin{equation}
\a(n) \a(n+2) \geq \a(n+1)^{2}. \notag
\end{equation}
\item[(C1)] There exists a constant $c_{1}$ such that for
all $n\leq m$,
\begin{equation}
\a(n) \leq c_{1}^{m} \a(m).  \notag
\end{equation}
\item[(C2)] There exists a constant $c_{2}$ such that for
all $n$ and $m$,
\begin{equation}
\a(n+m) \leq c_{2}^{n+m} \a(n) \a(m).  \notag
\end{equation}
\item[(C3)] There exists a constant $c_{3}$ such that for
all $n$ and $m$,
\begin{equation}
\a(n) \a(m) \leq c_{3}^{n+m} \a(n+m).  \notag
\end{equation}
\end{itemize}

Cochran et al.~\cite{cks} 
assumed condition (A1) with $\s=1$. But our (A1) is strong
enough to imply that the space of test functions is
contained in the $L^{2}$-space of the white noise measure.
In \cite{cks} conditions (A2) (B1) (B2) are considered.
Condition (A2) is to assure that the function $G_{\a}$
is an entire function.
Condition (B1) is used for the characterization theorem of
generalized functions in Theorems 5.1 and 6.1 \cite{cks}.
Condition (B2) is shown to imply condition (B1) in
Theorem 4.3 \cite{cks}.

In the papers by Asai et al.~\cite{akk1} \cite{akk2},
conditions ($\wt{\text{A}}2$) ($\wt{\text{B}}1$)
($\wt{\text{B}}2$) (B3) are considered. It can be
easily checked that condition (A1) implies condition
($\wt{\text{A}}2$). Condition ($\wt{\text{A}}2$) is to
assure that the function $G_{1/\a}$ 
is an entire function. In \cite{akk2}
condition ($\wt{\text{B}}1$) is used for the
characterization theorem of test functions. Condition
($\wt{\text{B}}2$) implies condition ($\wt{\text{B}}1$),
while obviously condition (B3) implies condition
($\wt{\text{B}}2$).

In the paper by Kubo et al.~\cite{kks}, conditions (C1)
(C2) (C3) are assumed in order to carry out the
distribution theory for a CKS-space. As pointed out in
\cite{kks}, condition (C3) implies condition (C1).

An important example of $\{\a(n)\}$ is the sequence
$\{b_{k}(n)\}$ of Bell's numbers of order $k\geq 2$. The
sequence $\{b_{k}(n)\}$ satisfies conditions (A1) (A2)
(B1) (as shown in \cite{cks}), (B2) (B3) (as shown in
\cite{akk1}) (C1) (C2) (C3) (as shown in \cite{kks}.)
Therefore, Bell's numbers satisfy all conditions in
the above list.

The essential conditions for distribution theory on a
CKS-space are (A1) (A2) (B2) ($\wt{\text{B}}2$) (C2)
(C3). All other conditions can be derived from these
six conditions except for (B3). We have taken
($\wt{\text{B}}2$) instead of (B3) for the following
reason. The condition (B3) is rather strong and we do
not know how to prove this condition for a growth
function $u$. Fortunately, we do
not need (B3) for white noise distribution theory.

\bigskip
\noindent
{\bf Acknowledgements.} N.~Asai is grateful for reserch support from 
the Daiko Foundation and the Kamiyama Foundation.
He also wants to thank Professor T. Sawada,
the Director of the  
International Institute of Advanced Studies,
for his constant 
encouragement.
H.-H.~Kuo is grateful for financial supports from the
Academic Frontier in Science (AFS) of Meijo University and
the Luso-American Foundation. He wants to thank AFS and
CCM, Universidade da Madeira for the warm hospitality during
his visits (February 15--21, 1998 and March 1--7, 1999 to
AFS, July 22-August 20, 1999 to CCM.) In particular, he gives
his deepest appreciation to Professors T. Hida and K. Sait\^o
(AFS) and M. de Faria and L. Streit (CCM) for arranging
the visits. H.-H.~Kuo also wants to thank Professor Y.-J.~Lee
for arranging his visit to Cheng Kung University in the
spring of 1998. At that time this joint research project
started. He thanks the financial support from the National
Science Council of Taiwan.


\end{document}